\documentclass[12pt,a4paper]{article}
\usepackage{amsmath,amssymb,amsfonts}
\def\Bbb{\mathbb}

\title{\bf Congruences modulo 4 for the alternating sum of the partial quotients}

\author{Kurt Girstmair\\ \small
Institut f\"ur Mathematik, Universit\"at Innsbruck   \\
\small Technikerstr. 13/7, A-6020 Innsbruck, Austria \\
\small Kurt.Girstmair@uibk.ac.at}

\date{}

\normalsize

\makeatletter
\let\@@maketitle=\maketitle
\def\maketitle{\def\thispagestyle##1{\relax}\@@maketitle}
\makeatother
\newtheorem{theorem}{Theorem}

\newtheorem{corollary}{Corollary}

\def\BE{\begin{equation}}
\def\EE{\end{equation}}
\def\BD{\begin{displaymath}}
\def\ED{\end{displaymath}}
\def\BA{\begin{array}}
\def\EA{\end{array}}
\def\BEA{\begin{eqnarray*}}
\def\EEA{\end{eqnarray*}}
\def\BI{\bibitem}

\def\Z{\Bbb Z}

\def\R{\Bbb R}

\def\phi{\varphi}

\def\MB{\mbox}
\def\LD{\ldots}

\def\BQ{``}
\def\EQ{'' }
\def\EQP{''}

\def\DED{Dedekind }

\begin{document}
\maketitle

\begin{abstract}
\noindent
We apply a technique used in $[$Tsukerman, Equality of Dedekind sums mod $\Z$, $2\Z$ and $4\Z$, ArXiv, Article-id 1408.3225]
combined with the Barkan-Hickerson-Knuth-formula in order to obtain congruences modulo $4$ for the alternating sum of the
partial quotients of the continued fraction expansion of $a/b$, where $0<a<b$ are integers.

\end{abstract}

\section*{Introduction and result}

Let $a$ and $b$ be positive integers, $a<b$. We consider the regular continued fraction expansion
\BD
 \frac ab=[0,a_1,\LD,a_n],
\ED
where all partial quotients $a_1,\LD,a_n$ are positive integers. We do not demand $a_n\ge 2$ but require $n$ to be odd instead.
So if $n$ is even and $a_n\ge 2$, we write $a/b=[0,a_1,\LD,a_n-1,1]$, and if $n$ is even and $a_n=1$, we write
$a/b=[0,a_1,\LD,a_{n-1}+1]$.

Let $a^*$ denote the inverse of $a$ mod $b$, i.e., $0<a<b$ and $aa^*\equiv 1$ mod $b$. In addition, let the positive integer $k$ be such that
\BD
   aa^*=1+kb.
\ED

In this note we express the congruence class of
\BD
   T(a,b)=\sum_{j=1}^n(-1)^{j-1}a_j
\ED
mod 4 and of
\BD
  D(a,b)=\sum_{j=1}^na_j
\ED
mod 2 in terms of $b$ and $k$. More precisely, we use the technique applied in \cite{Ts} in combination with the Barkan-Hickerson-Knuth-formula and obtain

\begin{theorem} 
\label{t1}
\hspace{1mm}

 {\rm (1)} If $a$ or $a^*$ is $\equiv 1$ mod $4$, then $T(a,b)\equiv b-k$ mod $4$.

 {\rm (2)} If $a$ or $a^*$ is $\equiv 3$ mod $4$, then $T(a,b)\equiv 2+k-b$ mod $4$.

 {\rm (3)} If $a$ or $a^*$ is $\equiv 2$ mod $4$, then $D(a,b)\equiv (b-k)/2$ mod $2$.

\end{theorem} 

Assertions (1) and (2) of Theorem \ref{t1} immediately yield

\begin{corollary} 
\label{c2}
If $a$ or $a^*$ is odd, then $D(a,b)\equiv b-k$ mod $2$.

\end{corollary} 

Of course, Theorem \ref{t1} and Corollary \ref{c2} do not cover all possible cases. For this reason we present the following

\medskip
\noindent
{\em Conjectures.} Suppose $a\equiv a^*\equiv 0$ mod 2.

   (1) If $a$ or $a^*$ is $\equiv 2$ mod $4$, then $T(a,b)\equiv (b-k)/2$ mod $4$.

   (2) If $a$ and $a^*$ are both $\equiv 0$ mod $4$, then $T(a,b)\equiv (k-b)/2$ mod $4$.

   (3) If $a$ and $a^*$ are both $\equiv 0$ mod $4$, then $D(a,b)$ is odd.

\medskip
\noindent
Again, Conjecture (3) is an immediate consequence of Conjecture (2). However, our method does not suffice to prove these conjectures.

\medskip
\noindent
{\em Remark.} Suppose we replace the above requirement \BQ $n$ odd\EQ by \BQ$a_n\ge 2$\EQP. Let $n$ be even and put
\BD
T'(a,b)=\sum_{j=1}^n(-1)^{j-1}a_j.
\ED
Then $T(a,b)=T'(a,b)+2$. Hence our results can easily be rephrased for this situation.

\section*{Proof of Theorem \ref{t1}}

Let the above notations hold.
The classical {\em \DED sum} $s(a,b)$ is defined by
\BE
\label{0.2}
   s(a,b)=\sum_{j=1}^b ((j/b))((aj/b))
\EE
where $((\LD))$ is the \BQ sawtooth function\EQ defined by
\BD
  ((t))=\left\{\begin{array}{ll}
                 t-\lfloor t\rfloor-1/2 & \MB{ if } t\in\R\smallsetminus \Z; \\
                 0 & \MB{ if } t\in \Z
               \end{array}\right.
\ED
(see, for instance, \cite[p. 1]{RaGr}). The Barkan-Hickerson-Knuth-formula says
\BD
 12s(a,b)=T(a,b)+\frac{a+a^*}b-3
\ED
(see \cite{Ba, Hi, Kn}; observe that $n$ is odd). On the other hand, the Dedekind sum is closely related to $I(a,b)$, the number of inversions of $a$ mod $b$, which
is defined in \cite{Ts}. Indeed, a result of C. Meyer says
\BD
   12bs(a,b)=-4I(a,b)+(b-1)(b-2)
\ED
(see \cite[Theorem 1.2]{Ts}).
Furthermore, it is an immediate consequence of a result of H. Sali\'e that the integer $I(a,b)$ satisfies
\BD
  4aI(a,b)\equiv (a-1)(b-1)(a+b-1) \MB{ mod }4b
\ED
(see \cite[Theorem 1.3]{Ts}). These three formulas give
\BD
  abT(a,b)\equiv -a^2-aa^*+3ab -(a-1)(b-1)(a+b-1)+a(b-1)(b-2) \MB{ mod } 4b.
\ED
If we write $aa^*=1+kb$ and expand the products on the right hand side, we obtain
\BD
  abT(a,b)\equiv -a^2b+3ab-2b+b^2-kb \MB{ mod } 4b
\ED
and
\BD
  aT(a,b)\equiv -a^2+3a-2+b-k \MB{ mod } 4.
\ED
If $a\equiv 1$ mod $4$, this gives
\BD
  T(a,b)\equiv b-k \MB{ mod } 4,
\ED
and in the case $a\equiv 3$ mod $4$
\BD
 -T(a,b)\equiv -2+b-k \MB{ mod 4}.
\ED
Now we observe $s(a,b)=s(a^*, b)$ (see \cite[p. 26]{RaGr}) and $(a^*)^*=a$. Thereby,
$T(a,b)=T(a^*,b)$. Altogether, we obtain assertions (1) and (2) of Theorem \ref{t1}.

If $a \equiv 2$ mod $4$, we have
\BD
   2T(a,b)\equiv b-k \MB{ mod }4
\ED
and
\BD
  T(a,b)\equiv D(a,b)\equiv (b-k)/2 \MB{ mod } 2.
\ED


\end{document}